\documentclass[11pt]{article}
\usepackage{amssymb}
\usepackage{amsbsy}
\usepackage{epsfig}
\usepackage{tikz}
\usepackage{amsmath}

\begin{document}
\pagestyle{myheadings}
\title{An explicit formula for the distance characteristic polynomial of threshold graphs }
\author{
   {\em Jo\~ao Lazzarin, Oscar F. M\'arquez and Fernando C. Tura}\thanks{
        Corresponding author. Email addresses: {\tt lazzarin@smail.ufsm.br} (J. Lazzarin), {\tt oscar.f.marquez-sosa@ufsm.br} (O.F. M\'arquez), {\tt ftura@smail.ufsm.br} (F.C. Tura).} \\
        Departamento de Matem\'atica\\
    Universidade Federal de Santa Maria, UFSM \\
    Santa Maria, RS, 97105-900, Brazil\\
}


\def\floor#1{\left\lfloor{#1}\right\rfloor}

\newcommand{\diagonal}[8]{
\begin{array}{| c | r |}
b_i & d_i \\
\hline
 #1 & #5 \\
        & \\
 #2 & #6 \\
        & \\
 #3  & #7 \\
         &  \\
 #4  & #8 \\
 \hline
 \end{array}
}

\newcommand{\thresholdmatrix}[8]{
\left [
     \begin{array}{ccccccc}
            &   &   &   &  #1   &   & #2 \\
            &   &   &   &  \vdots   &   & \vdots \\
            &   &   &   &  #1  &   & #2 \\
     \\
         #3  & \ldots   & #3 &   &  #4    &   & #5 \\
     \\
         #6  & \ldots   & #6  &   &  #7   &   &  #8 \\
    \end{array}
\right ]
}

\newcommand{\lambdamin}{\lambda_{\min, n}}
\newcommand{\formulamin}{
 \frac{1}{2} (  (
   \lfloor \frac{n}{3} \rfloor \! - \! 1 ) \! - \! \sqrt{ ( \lfloor \frac{n}{3} \rfloor \! - \! 1)^2 \! +\!  4
   (n \! - \! \lfloor \frac{n}{3} \rfloor )
   \lfloor \frac{n}{3} \rfloor
  } )
}
\newcommand{\casei}{{\bf Case~1}}
\newcommand{\caseii}{{\bf Case~2}}
\newcommand{\caseiii}{{\bf Case~3}}
\newcommand{\caseiv}{{\bf Case~4}}
\newcommand{\subia}{{\bf subcase~1a}}
\newcommand{\subib}{{\bf subcase~1b}}
\newcommand{\subic}{{\bf subcase~1c}}
\newcommand{\subiia}{{\bf subcase~2a}}
\newcommand{\subiib}{{\bf subcase~2b}}
\newcommand{\subiic}{{\bf subcase~2c}}
\newcommand{\subiiia}{{\bf subcase~3a}}
\newcommand{\subiiib}{{\bf subcase~3b}}
\newcommand{\subiiic}{{\bf subcase~3c}}
\newcommand{\subiva}{{\bf subcase~4a}}
\newcommand{\subivb}{{\bf subcase~4b}}
\newcommand{\subivc}{{\bf subcase~4c}}
\newcommand{\myvar}{x}
\newcommand{\exvar}{\frac{\sqrt{3} + 1}{2}}
\newcommand{\Prf}{{\bf Proof: }}
\newcommand{\PrfSketch}{{\bf Proof (Sketch): }}
\newcommand{\boldQ}{\mbox{\bf Q}}
\newcommand{\boldR}{\mbox{\bf R}}
\newcommand{\boldZ}{\mbox{\bf Z}}
\newcommand{\boldc}{\mbox{\bf c}}
\newcommand{\sign}{\mbox{sign}}
\newcommand{\alphaseq}{{\pmb \alpha}_{G,\myvar}}
\newcommand{\alphaseqGprime}{{\pmb \alpha}_{G^\prime,\myvar}}
\newcommand{\alphaseqlam}{{\pmb \alpha}_{G,-\lambdamin}}
\newtheorem{Thr}{Theorem}
\newtheorem{Pro}{Proposition}
\newtheorem{Que}{Question}
\newtheorem{Con}{Conjecture}
\newtheorem{Cor}{Corollary}
\newtheorem{Lem}{Lemma}
\newtheorem{Fac}{Fact}
\newtheorem{Ex}{Example}
\newtheorem{Def}{Definition}
\newtheorem{Prop}{Proposition}
\def\floor#1{\left\lfloor{#1}\right\rfloor}

\newenvironment{my_enumerate}{
\begin{enumerate}
  \setlength{\baselineskip}{14pt}
  \setlength{\parskip}{0pt}
  \setlength{\parsep}{0pt}}{\end{enumerate}
}
\newenvironment{my_description}{
\begin{description}
  \setlength{\baselineskip}{14pt}
  \setlength{\parskip}{0pt}
  \setlength{\parsep}{0pt}}{\end{description}
}

\maketitle

\begin{abstract}
 A threshold graph $G$ on $n$ vertices is defined by binary sequence of length $n.$
 In this paper we present an explicit formula for computing the distance characteristic polynomial 
 of a threshold graph from its binary sequence. As application, we show a several of  
 nonisomorphic  pairwise threshold graphs which are  $\mathcal{D}$-cospectral graphs. 
\end{abstract}
{\bf keywords:} threshold graph, distance characteristic polynomial, $\mathcal{D}$-cospectral graphs. 

\setlength{\baselineskip}{24pt}

\section{Introduction}
\label{intro}
All graphs considered in this paper are simple graphs, that is,  undirected, loop free and having no multiple edges.
Let  $G= (V,E)$ be a connected graph of order $n,$ where $V$  is the vertex set  and $E$ is the edge set.  The distance matrix $\mathcal{D}(G)$ of $G$ is an   $n\times n$ matrix $d_{ij}$ such that $d_{ij}$ is the distance (length of a shortest path) 
 between  $i$th and $j$th vertices in $G.$  

 The distance characteristic polynomial of $G,$ denoted by $P_{\mathcal{D}(G)}(x),$ can be expressed  as
$P_{\mathcal{D}(G)}(x) =\det(\mathcal{D}(G) - x I ).$ Their roots  are called distance eigenvalues of $G,$ or simply the $\mathcal{D}$-{\em spectrum} of $G.$ 
The distance eigenvalues were first studied by Graham and Pollack in 1971 to solve a data communication problem \cite{Graham,Indulal09}, and they have many applications on the literature \cite{Balaban,Dobrynin}.
We refer the reader to \cite{Aouchiche}  where several spectral  results were presented for
the distance matrix of graphs.

This paper is concerned  with threshold graphs, introduced by Chv\'atal and Hammer \cite{Chv} and Henderson and Zalcstein \cite{Hen} in 1977. They are an important class of graphs because of their numerous 
applications in diverse areas which include physics, biology and social sciences \cite{Mah95}. Threshold graphs can be characterized in many ways. One way of obtaining  a threshold graph is through a binary sequence that will be relevant to this paper, and we will describe in the next section.

There is a considerable body of knowledge on the spectral properties of threshold graphs related to adjacency matrix \cite{Bapat,Fritscher,furer,JTT2015,JTT2014,JTT2013,Stanic,SF2011,Simic}.
However, the literature does not seem to provide many results about the distance matrix of this class of graphs. One reason for this is due to the fact that distance matrix is dense while that adjacency matrix is relatively sparse. Thus the computation of the characteristic polynomial of distance matrix is computationally more complex problem.

 In this paper we attempt to fill this gap with presenting in Section \ref{Secformula} an explicit  formula for computing the distance characteristic polynomial of threshold graph from its binary sequence.
The distance eigenvalues $-2$ and $-1$ in threshold graphs  can be also obtained directly from its binary sequence. This is shown in Section \ref{Sec2}. 


 Two nonisomorphic graphs with the same spectrum are called {\em cospectral}. 
In recent years, there has been a growing interest to find families of cospectral graphs. There are many constructions in the literature \cite{godsil82,schwenk71}. This notion is originally defined for the adjacency matrix of the graph $G,$ but a natural extension of the problem is to find families of graphs that are cospectral  with relation to other matrices \cite{tura17}. As application, in Section \ref{Sec4},
 we show a several of  nonisomorphic  pairwise threshold graphs which are  $\mathcal{D}$-cospectral graphs.

\section{Preliminaries}
\label{Sec2}
 In this section, we present a formula for computing the multiplicities of distance eigenvalues $-2$ and $-1$ of a threshold graph,
as well as some known results.

\subsection{ The distance eigenvalues $-2$ and $-1.$} 
 Recall that  a vertex is  {\em isolated} if it has no neighbors, and is  {\em dominating} if it is  adjacent to all others vertices. 
 A threshold graph is obtained
through an iterative process which starts with an isolated vertex, and where, at each
step, either a new isolated vertex is added, or a {\em dominating} vertex (adjacent to all others vertices) is added.
 
 We represent a threshold graph $G$ on  $n$ vertices  using  a binary  sequence $(b_1, b_2, \ldots, b_n).$   Here $b_i=0$ if vertex $v_i$ was added as  isolated, and $b_i=1$  if vertex  $v_i$  was added  as a dominating vertex. We call our representation a {\em creation sequence}, and always take $b_1=0.$ If $n\geq 2,$  $G$ is connected if and only if $b_n=1.$

In constructing a distance matrix, we order the vertices
in the same way they are given in their creation sequence.
Figure~\ref{fig1} shows the
distance matrix $\mathcal{D}(G)$ of the
threshold graph $G$ represented by $(0,0,0,1,1,0,0,1)$ or $(0^3 1^2 0^2 1^1).$
\begin{figure}
 $$ 
 \small{
 \left[
  \begin{array}{cccccccc}
   0 & 2 & 2 & 1 & 1 &  2 & 2 & 1 \\
   2 & 0 & 2 & 1 & 1 &  2 & 2 & 1 \\
   2 & 2 & 0 & 1 & 1 &  2 & 2 & 1 \\
   1 & 1 & 1 & 0 & 1 &  2 & 2 & 1 \\
   1 & 1 & 1 & 1 & 0 &  2 & 2 & 1 \\
   2 & 2 & 2 & 2 & 2 &  0 & 2 & 1 \\
   2 & 2 & 2 & 2 & 2 &  2 & 0 & 1 \\
   1 & 1 & 1 & 1 & 1 &  1 & 1 & 0 \\
    \end{array}
  \right]
}$$
\caption{Distance matrix of threshold graph. \label{fig1} }
\end{figure}

 Let $m_{-2}(G)$ and $m_{-1}(G)$ denote the multiplicity of eigenvalues $-2$ and $-1,$ respectively,
in a threshold graph $G$ with distance matrix $\mathcal{D}.$
 We will represent a connected threshold graph by 
$G= (0^{a_1} 1^{a_2} 0^{a_3} \ldots 0^{a_{n-1}}1^{a_n})$ where each $a_i$ is a positive integer for $i = 1, \ldots,n.$

\begin{Lem}
\label{mainB} 
 For a connected threshold graph  $G= (0^{a_1} 1^{a_2} \ldots 0^{a_{n-1}} 1^{a_n})$
where each $a_i$ is a positive integer. Then   
 \begin{equation}
 \label{eq1}
 m_{-2}(G) = \sum_{i=1}^{\frac{n}{2}} (a_{2i-1} -1)
 \end{equation}
 and 
\begin{eqnarray}
\label{m(-1)}
 m_{-1}(G)=  \left\{\begin{array}{ccccc}
          \sum_{i=1}^{\frac{n}{2}} (a_{2i} - 1) & if&  a_1 >1 &  &  \\
                 &     &       &    \\
        1+  \sum_{i=1}^{\frac{n}{2}} (a_{2i} -1)      &    if    & a_1 =1  & & \\
   \end{array}\right.
 \end{eqnarray}  
\end{Lem}
\begin{Prf}
Let $G= (0^{a_1} 1^{a_2} \ldots 0^{a_{n-1}} 1^{a_n})$  be a connected threshold graph with distance matrix $\mathcal{D}(G).$ 
According Theorem 2.34  of \cite{Aouchiche}, if there are two vertices with the same neighborhood in a graph $G,$ then one root of distance polynomial is either $-1$ (if two vertices are adjacent) or $-2$ (if two vertices are not adjacent).

Follows that the distance eigenvalue $-2$ corresponds to set of vertices $0^{a_i},$ and hence its multiplicity is given by equation (\ref{eq1}). Similarly, the distance eigenvalue $-1$ corresponds to set of vertices $1^{a_i}$ and its multiplicity given by equation (\ref{m(-1)}). 
\end{Prf}

\subsection{ The parameter $\gamma_{n,l}^a$}

In order to obtain an explicit formula to $P_{\mathcal{D}(G)}(x)$  for $G= (0^{a_1} 1^{a_2} \ldots 0^{a_{n-1}} 1^{a_n})$ 
we need to introduce the parameter $\gamma_{n,l}^a.$

Let $[n]=\{1,2,\ldots n\}$, and let $I_{n,l}$ the set of increasing sequences   in $[n]$  of length $l$ such that if $(t_1,t_2,\ldots,t_l)\in I_{n,l}$ then $t_i\equiv n+i-l \, (mod\hspace{0,2cm} 2).$ In other words, elements in $I_{n,l}$ are increasing sequences alternating even and odds numbers such that the last term has the same than parity $n$. For instance
$$I_{7,4}=\{(2,3,4,5),(2,3,4,7),(2,3,6,7), (2,5,6,7),(4,5,6,7) \},$$ while 
$$I_{6,4}=\{(1,2,3,4),(1,2,3,6),(1,2,5,6),(1,4,5,6),(3,4,5,6) \}.$$
In general, for any $l$ the sequences in $I_{7,l}$ must finish in an odd number, while all the sequences in $I_{6,l}$ must finish in an even number.  Given a sequence $\mathbf{t}=(t_1,t_2,\ldots,t_l)$ we denote  $a_\mathbf{t}=a_{t_1} a_{t_2}\cdots a_{t_l}$. Based in this notation we have

\begin{Def}
Let $a = (a_1, a_2, \ldots, a_n) $  be fixed sequence of positive integers. We define the following parameter 

  $$   \gamma _{n,l}^a=  \left\{\begin{array}{ccc}
           \sum\limits_{\mathbf{t}\in I_{n,l}}  a_\mathbf{t}& if &1\leq l \leq n \\
                1  & if &    l=0 .     \\
   \end{array}\right.$$
 \end{Def}

We write $n=2 m+ r_0$, with $r_0\in \{0,1\}$ and  define $r_1\in\{0,1\}$ such that $r_1\equiv r_0+1 (mod \hspace{0,2cm}2)$. For $n\ge2$ and $y=x+1,$ 
the following result was given in \cite{JOF}.

\begin{Prop}\label{formula1}
Let $a=(a_1, a_2, \ldots, a_n)$  be fixed sequence of positive integers.The determinant of the following $n\times n$ tridiagonal matrix 
 \begin{eqnarray}
\label{eq11}
  \left[\begin{array}{cccccc}
                x+a_1   &   -y           &   0       &  0                  & 0 &  0  \\
                x   &  - a_2   &   -x   &      0                     & 0& 0\\
                0         &      y        &   a_3   &     -y           &    0  & 0  \\
               0 &  0         &    x &  -a_4                   & -x  & 0    \\
                \vdots     &         \ddots &        \ddots        & \ddots        &\ddots & \vdots     \\
              0              &           0     &         \ldots            &  y&  a_{n-1} & -y\\
              0       &    0     &     0  &  \ldots  &     x  &     -a_n \\
\end{array}\right]
\end{eqnarray}
can be viewed as a bivariate polynomial $p_n^{(a)}(x,y)$ with integer coefficients and computed by
\begin{align}\label{eq12}
	p_{n}^{(a)}(x,y)&= x^{r_0} \sum_{k=0}^{m} (-1)^{m-k}x^k y^k \gamma_{n, n-2 k-r_0}^a  +\\ \notag
	&   x^{r_1} \sum_{k=0}^{m-r_1} (-1)^{m-k}x^k y^k \gamma_{n, n-2k-r_1}^a.  \\ \notag
	\end{align}
\end{Prop}

\section{An explicit formula for $P_{\mathcal{D}(G)}(x)$ }
\label{Secformula}

In this section we present an explicit formula for the distance characteristic polynomial of threshold graph from its binary sequence.
A similar formula related to adjacency matrix was given in \cite{JOF}.

We begin  with an auxiliary result that will play an important role in the sequel.

\begin{Lem} 
\label{lem3}  Let $G= (0^{a_1} 1^{a_2} \ldots 0^{a_{n-1}} 1^{a_n})$ be a connected threshold graph with distance matrix $\mathcal{D}(G).$
Let $m_{-2}(G)$ and $m_{-1}(G)$ be the multiplicities of eigenvalues $-2$ and $-1$ of $G,$ respectively. The distance characteristic polynomial of $G$ is, to within a sign,
\begin{eqnarray}
\label{eq8}
   \begin{array}{c}
           P_{\mathcal{D}(G)}(x) = (x+2)^{m_{-2}(G)} (x+1)^{m_{-1}(G)} \mathcal{Q}(x),  \\
   \end{array}
 \end{eqnarray}  
  where $\mathcal{Q}(x)$  is the characteristic polynomial of the matrix below 

\begin{eqnarray}
\label{eq7}
 \mathcal{Q}= \left[\begin{array}{ccccccc}
                2(a_1 -1)   &   a_2            &   2 a_3       &  a_4               & \ldots    &2 a_{n-1}  &  a_{n}  \\
               a_1    &   a_2-1    &   2a_3   &      a_4            & \ldots           & 2 a_{n-1}& a_{n}\\
               2 a_1         &      2 a_2        &   2(a_3 -1)  &     a_4        & \ldots      &     2 a_{n-1}   &  a_n  \\
               a_1 &  a_2            &    a_3  &  a_4 -1              & \ldots      & 2 a_{n-1}   & a_{n}    \\
            \vdots         &         \vdots &        \vdots        & \vdots   &  \ddots       &\vdots & \vdots     \\
              2a_1              &           2a_2     &         2a_3      &     \ldots              &   & 2 ( a_{n-1} -1)& a_n\\
        a_{1}       &    a_{2}     &a_{3}  &\ldots & a_{n-2}  &  a_{n-1}   &     a_n -1\\
\end{array}\right]
\end{eqnarray}
\end{Lem}
\begin{Prf}
Let  $G= (0^{a_1} 1^{a_2} \ldots 0^{a_{n-1}} 1^{a_n})$ be a connected threshold graph with distance matrix $\mathcal{D}(G).$ Since that the multiplicities $m_{-2}(G)$ and $m_{-1}(G)$ can  be obtained from sequence binary of $G,$ according Lemma \ref{mainB}, we will determine the $\mathcal{Q}(x).$ Let $\lambda$ be an eigenvalue of $\mathcal{D}(G)$ with an eigenvector associated  $x=(x_1, x_2, \ldots, x_n).$  Using the system $\mathcal{D}x=\lambda x$ we obtain the following equations:
\begin{eqnarray}
\sum_{i=2}^{a_1}2x_{i} +\sum_{i=a_1 +1}^{a_2}x_{i} +\sum_{i=a_2 +1}^{a_3}2x_{i} + \ldots +\sum_{i=a_{n-1}+1}^{a_n}x_{i} =\lambda x_{1}\nonumber \\
2x_1 + \sum_{i=3}^{a_1}2x_{i} +\sum_{i=a_1 +1}^{a_2}x_{i} +\sum_{i=a_2 +1}^{a_3}2x_{i} + \ldots +\sum_{i=a_{n-1}+1}^{a_n}x_{i} =\lambda x_{2}\nonumber \nonumber\\
\vdots \nonumber \\
\sum_{i=1}^{a_1}x_{i} + \sum_{i=a_1 +1}^{a_2}x_{i}+ \sum_{i=a_2 +1}^{a_3}x_{i} + \ldots + \sum_{i=a_{n-1}+1 }^{a_n}x_{i}=\lambda x_{n} \nonumber
\end{eqnarray}
Subtracting the first two equations in this system,  we obtain that
$$-2x_1 +2x_2 = \lambda(x_1-x_2)$$
which is equivalent to
$$ (\lambda+2)(-x_1 +x_2)=0.$$
Assuming that $\lambda \neq -2,$ we have that $x_1=x_2.$  Using a similar argument for the first $a_1$ equations, we obtain that $x_3=x_1, x_4=x_1$ and so on. Then 
$$x=(\underbrace{x_1, x_1, \ldots,x_1}_{a_1 }, \underbrace{x'_2, x'_2 \ldots, x'_2}_{a_2},\underbrace{x'_3, x'_3 \ldots, x'_3}_{a_3},\ldots,a_1 x_1+a_2x'_2+ \ldots +(a_n-1) x_{n}).$$
Replacing this solution in the original system $\mathcal{D}x =\lambda x,$ we obtain the desire matrix (\ref{eq7}) and the result follows.
\end{Prf}\\
\noindent{\bf Remark:} An alternative approach to prove the previous Lemma is to use an equitable partition of vertex set  $V(G)$ and the divisor technique, according \cite{Atik}.

\begin{Lem}
\label{lem4} 
Let $n\geq 1$ be a positive integer and let
\label{eq11}
$$ M_n= \left[\begin{array}{ccccccc}
                -1   &   1           &   0       &  0                  & 0 &   & 0  \\
                -1   &  0   &   1  &      0                     & 0&      &0\\
                0         &      -1       &   0   &     1         &    0  &     &0  \\
                \vdots     &         \ddots &        \ddots        & \ddots        &\ddots &     &\vdots     \\
              0              &           0     &         \ldots            &  -1&  0 &   &1\\
              0       &    0     &     0  &  \ldots  &     -1 &     &0 \\
\end{array}\right]$$
be a tridiagonal matrix of order $n.$ Then the $det M_n = (-1)^{n}.$
\end{Lem}
\begin{Prf} We prove the result by induction on $n.$ The cases $n=1$ and $n=2$  is easy to verify. Assume the result to be true for $M_k,  2 \leq k \leq n-1.$ A simple Laplace expansion shows that
$$ det M_n = (-1)(-1)^{2n-1} \cdot det M_{n-2}$$
$$ = (-1)\cdot(-1)\cdot  (-1)^{n-2}= (-1)^n.$$
\end{Prf}

\begin{Thr}
\label{formula}
Let $G= (0^{a_1} 1^{a_2} \ldots 0^{a_{n-1}} 1^{a_n})$ be a connected threshold graph with distance matrix $\mathcal{D}(G).$
Let $m_{-2}(G)$ and $m_{-1}(G)$ be the multiplicities of eigenvalues $-2$ and $-1$ of $G,$ respectively. The distance characteristic polynomial of $G$ is, to within a sign, 
$$P_{\mathcal{D}(G)}(x) = (x+2)^{m_{-2}(G)} (x+1)^{m_{-1}(G)} \mathcal{Q}(x), \hspace{0,5cm} where $$ 
\begin{equation}
\label{eq13}
\mathcal{Q}(x)=- p_n^{(-a)}(z,y) +2y p_{n-1}^{(-a)}(z,y)
\end{equation}
with $p_n^{(-a)}(z,y)$  given by (\ref{eq12}), replacing $z=x+2, y=x+1$ and  each $a_i$ by $-a_i.$
\end{Thr}
\begin{Prf}
According Lemma \ref{lem3}, it is sufficient to show that $\mathcal{Q}(x)$ satisfies the equation (\ref{eq13}).
Let $\mathcal{Q}-xI$ be the matrix  where $\mathcal{Q}$ is the matrix given by (\ref{eq7}),
\begin{center}
$ \left[\begin{array}{cccccc}
               2(a_1-1) -x   &   a_2          &   2a_3                     & \ldots  & 2a_{n-1} &a_n  \\
               a_1    &   a_2-1-x    &   2a_3                 & \ldots         &  2a_{n-1}  & a_n\\
               2a_1         &      2a_2       &  2(a_3 -1) -x        & \ldots     &   2a_{n-1}  &     a_n \\
            \vdots         &         \vdots &                 \ddots &     & \vdots   & \vdots \\
              2a_1              &           2a_2     &         2a_3      &     \ldots           & 2(a_{n-1} -1)-x    & a_n\\
               a_1       &    a_2     &   a_3  &\ldots    &     a_{n-1}  &   a_n -1-x \\
\end{array}\right]$
\end{center}

and replacing $y=x+1$ and $z = x+2,$  giving the  matrix:
\begin{eqnarray}
\label{eq}
 \left[\begin{array}{cccccc}
               2a_1-z   &   a_2          &   2a_3                     & \ldots  & 2a_{n-1} &a_n  \\
               a_1    &   a_2-y    &   2a_3                 & \ldots         &  2a_{n-1}  & a_n\\
               2a_1         &      2a_2       &  2a_3 -z        & \ldots     &   2a_{n-1}  &     a_n \\
            \vdots         &         \vdots &                 \ddots &     & \vdots   & \vdots \\
              2a_1              &           2a_2     &         2a_3      &     \ldots           & 2a_{n-1} -z    & a_n\\
               a_1       &    a_2     &   a_3  &\ldots    &     a_{n-1}  &   a_n -y \\
\end{array}\right]
\end{eqnarray}

Let $M_n$ be the matrix of Lemma \ref{lem4}, and using that $det M_n = (-1)^n, $ follows

$M_n  \times \left[ 
\begin{array}{cccccc}
    2a_1-z   &   a_2          &   2a_3                     & \ldots  & 2a_{n-1} &a_n  \\
               a_1    &   a_2-y    &   2a_3                 & \ldots         &  2a_{n-1}  & a_n\\
               2a_1         &      2a_2       &  2a_3 -z        & \ldots     &   2a_{n-1}  &     a_n \\
            \vdots         &         \vdots &                 \ddots &     & \vdots   & \vdots \\
              2a_1              &           2a_2     &         2a_3      &     \ldots           & 2a_{n-1} -z    & a_n\\
               a_1       &    a_2     &   a_3  &\ldots    &     a_{n-1}  &   a_n -y \\\end{array}%
\right] $

\begin{eqnarray}
\label{eq10}
= \left[\begin{array}{ccccccc}
               z -a_1   &   -y        &  0                    & \ldots  &0 &0 & 0  \\
                z  &   a_2    &   -z                 & \ldots        &0 &  0  & 0\\
               0         &      y      &  -a_3         & \ldots    & 0 &   0 &     0\\
            \vdots         &         \vdots &                 \ddots &   &  & \vdots   & \vdots \\
              0            &           0    &         0    &     \ldots  &     y    & -a_{n-1}     & -y\\
               -2a_1       &    -2a_2     &   -2a_3  &\ldots    &-2a_{n-2} &    z- 2a_{n-1}  &   -a_n  \\
\end{array}\right]
\end{eqnarray}

Since that the determinant of matrix (\ref{eq10}) can be computed as determinant of following matrices

\begin{eqnarray}
\label{eqc}
= \left[\begin{array}{ccccccc}
               z -a_1   &   -y        &  0                    & \ldots  &0 &0 & 0  \\
                z  &   a_2    &   -z                 & \ldots        &0 &  0  & 0\\
               0         &      y      &  -a_3         & \ldots    & 0 &   0 &     0\\
            \vdots         &         \vdots &                 \ddots &   &  & \vdots   & \vdots \\
              0            &           0    &         0    &     \ldots  &     y    & -a_{n-1}     & -y\\
               -2a_1       &    -2a_2     &   -2a_3  &\ldots    &-2a_{n-2} &    2z- 2a_{n-1}  &  0 \\
\end{array}\right]
\end{eqnarray}
\begin{eqnarray}
\label{eqd}
- \left[\begin{array}{ccccccc}
               z -a_1   &   -y        &  0                    & \ldots  &0 &0 & 0  \\
                z  &   a_2    &   -z                 & \ldots        &0 &  0  & 0\\
               0         &      y      &  -a_3         & \ldots    & 0 &   0 &     0\\
            \vdots         &         \vdots &                 \ddots &   &  & \vdots   & \vdots \\
              0            &           0    &         0    &     \ldots  &     y    & -a_{n-1}     & -y\\
               0       &    0    &   0  &\ldots    &0 &    z  &  a_n \\
\end{array}\right]
\end{eqnarray}
and by Proposition \ref{formula1}, the determinant of matrix (\ref{eqd})  is $-p_{n}^{(-a)}(z,y).$
Now, we show that the determinant of matrix (\ref{eqc})  is $2yp_{n-1}^{(-a)}(z,y).$
Using properties of determinant, we have that 
\begin{eqnarray}
\label{eqe}
= -2 \left[\begin{array}{ccccccc}
               z -a_1   &   -y        &  0                    & \ldots  &0 &0 & 0  \\
                z  &   a_2    &   -z                 & \ldots        &0 &  0  & 0\\
               0         &      y      &  -a_3         & \ldots    & 0 &   0 &     0\\
            \vdots         &         \vdots &                 \ddots &   &  & \vdots   & \vdots \\
              0            &           0    &         0    &     \ldots  &     y    & -a_{n-1}     & -y\\
               a_1       &    a_2     &   a_3  &\ldots    &a_{n-2} &    a_{n-1}-z  &  0 \\
\end{array}\right]
\end{eqnarray}
and by Laplace expansion in the $n$th column of matrix (\ref{eqe})
\begin{eqnarray}
\label{eqf}
= -2 (-y)(-1)^{2n-1}\left[\begin{array}{cccccc}
               z -a_1   &   -y        &  0                    & \ldots  &0 &0   \\
                z  &   a_2    &   -z                 & \ldots        &0 &  0  \\
               0         &      y      &  -a_3         & \ldots    & 0 &   0 \\
            \vdots         &         \vdots &                 \ddots &   &  \vdots & \vdots    \\
                 0           &       0   &  \ldots     &    z    &  a_{n-2}   &  -z \\
               a_1       &    a_2     &   a_3  &\ldots    &a_{n-2} &    a_{n-1}-z   \\
\end{array}\right]
\end{eqnarray}
Performing the following operations $R_{n-1} \leftarrow \sum_{i=1}^{n-2} (-1)^{i+1}R_i +R_{n-1},$ giving 
\begin{eqnarray}
\label{eqf}
= 2y\left[\begin{array}{cccccc}
               z -a_1   &   -y        &  0                    & \ldots  &0 &0   \\
                z  &   a_2    &   -z                 & \ldots        &0 &  0  \\
               0         &      y      &  -a_3         & \ldots    & 0 &   0 \\
            \vdots         &         \vdots &                 \ddots &   &  \vdots & \vdots    \\
                 0           &       0   &  \ldots     &    z    &  a_{n-2}   &  -z \\
               0       &    0    &   0  &\ldots    & y &    -a_{n-1}  \\
\end{array}\right]
\end{eqnarray}
which the determinant is  $2yp_{n-1}^{(-a)}(z,y),$
and the result follows.
\end{Prf}

\begin{Ex}
We apply the formula given in Theorem \ref{formula} to the threshold graph $G=(0^{a_1} 1^{a_2} 0^{a_3} 1^{a_4})$ with $a_1>1.$ According Theorem \ref{mainB} the multiplicities of $-2$ and $-1$ are
$m_{-2}(G) = \sum_{i=1}^{2} (a_{2i-1} -1)$ and $m_{-1}(G) = \sum_{i=1}^{2} (a_{2i} -1).$ 
Since $P_{\mathcal{D}(G)}(x) = (x+2)^{m_{-2}(G)} (x+1)^{m_{-1}(G)} \mathcal{Q}(x)$ where $\mathcal{Q}(x)= -p_4^{(-a)}(z,y) +2yp_3^{(-a)}(z,y),$ 
by Proposition \ref{formula1}, follows
\begin{align*}
p_4^{(a)}(z,y) &=z^{0} \sum\limits_{k=0}^{2} (-1)^{1-k}z^k y^k \gamma_{4, 4-2k-0}^a  +\,  z^{1} \sum\limits_{k=0}^{2-1} (-1)^{1-k}z^k y^k \gamma_{4,4-2k-1}^a\\
&=  z^0 \left( -\gamma_{4 ,4}^a +z y \gamma_{4 ,2}^a - z^2 y^2 \gamma_{4,0}^a  \right)  +
  z^{1} \big(-\gamma_{4,3}^a+z y \gamma_{4,1}^a \big)\\
&= z^2 y^2 -(a_2 +a_4) z^2 y -zy( a_1a_2 + a_1 a_4 +a_3 a_4)  +z(a_2 a_3 a_4) -a_1a_2 a_3 a_4.
\end{align*}
and
\begin{align*}
p_3^{(a)}(z,y) &=z^{1} \sum\limits_{k=0}^{1} (-1)^{1-k}z^k y^k \gamma_{3,3-2 k-1}^a  +\,  z^{0} \sum\limits_{k=0}^{1} (-1)^{1-k}z^k y^k \gamma_{3,3-2k}^a\\
&=  z^1 \left( -\gamma_{3,2}^a +z y \gamma_{3,0}^a   \right)  +  z^{0} \big(-\gamma_{3,3}^a+zy \gamma_{3,1}^a \big)\\
&= z^2y +zy (a_1 +a_3) -z(a_2a_3) -a_1 a_2 a_3.
\end{align*}
\noindent Replacing $z=x+2,y=x+1$ and each $a_i$ by $-a_i,$
then $P_{\mathcal{D}(G)}(x)$ is 
\begin{align*}
P_{\mathcal{D}(G)}(x) &=  (x+2)^{m_{-2}(G)} (x+1)^{m_{-1}(G)}   \{-x^4 +x^3(2a_{1}+a_{2}+2a_{3}+a_{4}-6)\\
& +x^2(8a_{1}+5a_{2}+8a_{3}+5a_{4}-a_{1}a_{2}-a_{1}a_{4}+2a_{2}a_{3}-a_{3}a_{4}-13) \\
&+x(10a_{1}+8a_{2}+10a_{3}+8a_{4}-3a_{1}a_{2}-3a_{1}a_{4}+6a_{2}a_{3}-3a_{3}a_{4}-2a_{1}a_{2}a_{3}\\
&-a_{2}a_{3}a_{4}-12)   + 4a_{1}+4a_{2}+4a_{3}+4a_{4}-2a_{1}a_{2}-2a_{1}a_{4}+4a_{2}a_{3}-2a_{3}a_{4} \\
& -2a_{1}a_{2}a_{3}-2a_{2}a_{3}a_{4}+a_{1}a_{2}a_{3}a_{4}-4 \}.
\end{align*}
\end{Ex}

\section{$\mathcal{D}$-cospectral graphs}
\label{Sec4}

In this section, we present some connected threshold graphs which are $\mathcal{D}$-cospectral graphs.

\begin{Lem}
\label{Lema bom}
If $G=(0^{a_{1}}1^{a_{2}}0^{a_{3}}\ldots 0^{a_{n-1}}1^{a_{n}}) 
$ and $G^{\prime
}=(0^{b_{1}}1^{b_{2}}0^{b_{3}}\ldots 0^{b_{n-1}}1^{b_{n}}) 
$ are  $\mathcal{D}$-cospectral  then 
\begin{eqnarray*}
\gamma _{n,l}^{a} &=&\gamma _{n,l}^{b} \\
\gamma _{n-1,l}^{a} &=&\gamma _{n-1,l}^{b}
\end{eqnarray*}%
for $l=1,2.$
\end{Lem}

\begin{Prf}
We note that  $\gamma _{n,1}^{a}=m_{-2}(G)=\gamma _{n,1}^{b}$ and $\gamma
_{n-1,1}^{a}=m_{-1}(G)=\gamma _{n-1,1}^{b}$ and since that the coefficient of $x^{n-2}$
in $Q^{a}(y+1,y)=Q^{b}(y+1,y)$ are equal  if and only if $-\gamma
_{n,2}^{a}+2\gamma _{n-1,2}^{a}=-\gamma _{n,2}^{b}+2\gamma _{n-1,2}^{b}$. \
By other hand we have that  $\gamma _{n,2}^{a}+\gamma _{n-1,2}^{a}=\gamma
_{n,1}^{a}.\gamma _{n-1,1}^{a}$, and the result follows.
\end{Prf}

\begin{Lem}
\label{Lema a1}
If $G=(0^{a_{1}}1^{a_{2}}0^{a_{3}}1^{a_{4}}) 
$ and $G^{\prime
}=(0^{b_{1}}1^{b_{2}}0^{b_{3}}1^{b_{4}}) 
$ are  $\mathcal{D}$-cospectral  then 
\begin{equation}
\label{eq16}
a_{1}a_{4}+2a_{1} =b_{1}b_{4}+2b_{1}
\end{equation}
\begin{equation}
\label{eq17}
2a_{1}+a_{4} =2b_{1}+b_{4}
\end{equation}

\end{Lem}

\begin{Prf}
For $G=(0^{a_{1}}1^{a_{2}}0^{a_{3}}1^{a_{4}})$ and $G^{^{\prime
}}=(0^{b_{1}}1^{b_{2}}0^{b_{3}}1^{b_{4}})$ to be cospectral, we should have  
$P_{\mathcal{D}(G)}^{a}(x)=P_{\mathcal{D}(G')}^{b}(x).$ From example 1, follows  (taking $x=y-1)$

$Q^{(a_{1},a_{2},a_{3},a_{4})}(y)=\allowbreak \left(
2a_{1}+a_{2}+2a_{3}+a_{4}-2\right) y^{3}-y^{4}+\left(
2a_{1}+2a_{2}+2a_{3}+2a_{4}\right. $\\
$-a_{1}a_{2}-a_{1}a_{4}+2a_{2}a_{3}-a_{3}a_{4}-1) y^{2}$
$+\left(
a_{2}+a_{4}-a_{1}a_{2}-a_{1}a_{4}+2a_{2}a_{3}-a_{3}a_{4}-2a_{1}a_{2}a_{3}-a_{2}a_{3}a_{4}\right) \allowbreak y+\left( a_{1}a_{2}a_{3}a_{4}-a_{2}a_{3}a_{4}\right).
$

By Lemma \ref{Lema bom}, if $\mathcal{Q}^{a}(y)=\mathcal{Q}^{b}(y)$ then

\begin{eqnarray*}
a_{1}a_{2}a_{3}a_{4}-a_{2}a_{3}a_{4} &=&b_{1}b_{2}b_{3}b_{4}-b_{2}b_{3}b_{4}
\\
2a_{1}a_{2}a_{3}+a_{2}a_{3}a_{4} &=&2b_{1}b_{2}b_{3}+b_{2}b_{3}b_{4}
\end{eqnarray*}
adding this  two  equations, we obtain
\begin{eqnarray*}
(a_{1}a_{4}+2a_{1})a_{2}a_{3} &=&(b_{1}b_{4}+2b_{1})b_{2}b_{3} \\
(2a_{1}+a_{4})a_{2}a_{3} &=&(2b_{1}+b_{4})b_{2}b_{3}.
\end{eqnarray*}
Since that  $\gamma _{3,2}^{a}=$\bigskip $\gamma
_{3,2}^{b}=a_{2}a_{3}=b_{2}b_{3},$ 
 a simple division lead us to the result.
\end{Prf}

\begin{Thr}
$G=(0^{a_{1}}1^{a_{2}}0^{a_{3}}1^{a_{4}})$ and $G^{^{\prime
}}=(0^{b_{1}}1^{b_{2}}0^{b_{3}}1^{b_{4}})$ are nonisomorphic and $\mathcal{D}$-cospectral graphs if and only if
the following holds:
\begin{enumerate}
\item $\alpha = a_{1} - b_{1} > 0$ and $a_{1}\neq 1$ and  $b_{1}\neq 1;$ 

\item $b_{2}$  is even;

\item $b_{3}=\frac{b_{2}}{2}+\alpha ;$

\item $b_{4}=2(a_{1}-1);$

\item $a_{2}=b_{2}+2\alpha ;$

\item $a_{3}=\frac{b_{2}}{2};$

\item $a_{4}=2(b_{1}-1).$
\end{enumerate}
\end{Thr}

\begin{Prf}
We assume that $G=(0^{a_{1}}1^{a_{2}}0^{a_{3}}1^{a_{4}})$ and $G^{^{\prime
}}=(0^{b_{1}}1^{b_{2}}0^{b_{3}}1^{b_{4}})$ are nonisomorphic and $\mathcal{D}$-cospectral graphs.
By Lemma \ref{Lema a1},  $a_{1}=1$ if and only if $b_{1}=1.$ Now, $a_{1}=b_{1},$ again using the Lemma \ref{Lema a1} follows that $a_{4}=b_{4}$
and consequently $a_{2}=b_{2}$ and $%
b_{3}=a_{3},$ contradicting   $G\neq G^{\prime }$.


For the itens (iv) and (vii), by Lemma \ref{Lema a1} and Equation (\ref{eq17}), $a_4 = b_4 -\alpha,$ for $\alpha = a_1 - b_1.$ Replacing this in Equation (\ref{eq16})  we have that  $(a_{1}\left( b_{4}-2\alpha \right)
+2a_{1})=(b_{1}b_{4}+2b_{1})\Rightarrow \alpha b_{4}-2\alpha a_{1}+2\alpha
=0\Rightarrow \alpha =0$ or $b_{4}=2a_{1}-2$ and $a_{4}=2a_{1}+2+2\alpha
=2b_{1}-2.$ 

For the remainder items, we fixe $\alpha =a_{1}-b_{1}>0.$ According Lema \ref{Lema bom} we have that   $b_{3}-a_{3}=\alpha $ and  $a_{3}a_{2}=b_{3}b_{2},$
and by Equation (\ref{eq17}) results $a_{2}-b_{2}=2\alpha.$ 
Using similar procedures we obtain the other statements.

Now we assume that the items (i), (ii), $\ldots$ (vii) hold. We will prove
$G=(0^{a_{1}}1^{a_{2}}0^{a_{3}}1^{a_{4}})$ and $G^{^{\prime
}}=(0^{b_{1}}1^{b_{2}}0^{b_{3}}1^{b_{4}})$ are $\mathcal{D}$-cospectral graphs.

It is clear that  $m_{-1}(G) = a_{2}+a_{4} -2 =b_{2}+2(a_{1} -b_1) +2(b_1-1) -2 =b_{2}+b_{4}- 2= m_{-1}(G^{^{\prime}}).$
Similarly, we obtain that $m_{-2}(G) =  m_{-2}(G^{^{\prime}}).$

 For computing the remainder terms, we fixe  
$b_{2}=2\beta,$ for $\beta >0$, and $\alpha >0,$ then
$a=(\alpha +b_{1},2(\alpha +\beta ),\beta ,2(b_{1}-1)),$ and 
$b=(b_{1},2\beta ,\alpha +\beta ,2\alpha +2b_{1}-2)$, which provide

$\mathcal{Q}^{a}(x)=( 4\alpha +4\beta +4b_{1}-8)
x^{3}-x^{4}+\left( 2\alpha \beta -2\alpha ^{2}-4\alpha b_{1}+20\alpha
+4\beta ^{2}-4\beta b_{1}+20\beta \right. $
$ -2b_{1}^{2}+20b_{1}-23)
x^{2} + ( 10\alpha \beta -6\alpha ^{2}-4\alpha \beta ^{2}-8\alpha \beta
b_{1}-4\alpha ^{2}\beta -12\alpha b_{1}+32\alpha -8\beta ^{2}b_{1}$
$+16\beta ^{2}$
$-12\beta b_{1}+32\beta -6b_{1}^{2}+32b_{1}-28)x
+4\alpha ^{2}\beta b_{1}-8\alpha ^{2}\beta -4\alpha ^{2}+4\alpha \beta
^{2}b_{1}$
$-8\alpha \beta ^{2}+4\alpha \beta b_{1}^{2}-16\alpha \beta
b_{1}+12\alpha \beta -8\alpha b_{1}+16\alpha +4\beta ^{2}b_{1}^{2}-16\beta
^{2}b_{1}+16\beta ^{2}-8\beta b_{1}+16\beta -4b_{1}^{2}+16b_{1}-12=\mathcal{Q}^{b}(x),$
and the results follows.
\end{Prf}

\begin{Cor}
\label{Cor1} For positive integers $i, j, k$ and $l,$ the connected 
threshold graphs $G=(0^i 1^{2j} 0 ^k 1^{2l})$ and $G^{^{\prime
}}=(0^{l+1} 1^{2k} 0^{j} 1^{2i-1})$  are nonisomorphic and $\mathcal{D}$-cospectral graphs, if $i+k=l+j, i > 1,$ and $ l >1.$
\end{Cor}


\section{Acknowledgment}

F.C. Tura acknowledges the support of FAPERGS (Grant 17/2551-0000813-8).

\setlength{\baselineskip}{14pt}

\end{document}